\documentclass{article}

\usepackage{graphicx}
\usepackage{amsmath}

\newtheorem{theorem}{Theorem}

\newtheorem{corollary}[theorem]{Corollary}

\newtheorem{definition}[theorem]{Definition}

\newtheorem{lemma}[theorem]{Lemma}

\newtheorem{proposition}[theorem]{Proposition}
\newtheorem{remark}[theorem]{Remark}

\newenvironment{proof}[1][Proof]{\textbf{#1.} }{\ \rule{0.5em}{0.5em}}

\begin{document}

\title{Morphisms which are continuous on a neighborhood of the base of a groupoid%
\thanks{%
This work is supported by MEdC-CNCSIS grant At 127/2004.}}
\author{M\u{a}d\u{a}lina Roxana Buneci \\
University Constantin Br\^{a}ncu\c{s}i of T\^{a}rgu-Jiu}
\date{}
\maketitle

\begin{abstract}
Kirill Mackenzie raised in \cite{ma} (p. 31) the following question: given a
morphism $F:\Omega \rightarrow \Omega ^{\prime }$, where $\Omega $ and $%
\Omega ^{\prime }$ are topological groupoids and $F$ is continuous on a
neighborhood of \ the base in $\Omega $, is it true that $F$ is
continuous everywhere?

This paper gives a negative answer to that question. Moreover, we shall
prove that for a locally compact groupoid $\Omega $ with non-singleton
orbits and having open target projection, if we assume that the continuity
of every morphism $F$ on the neighborhood of \ the base in $\Omega $
implies the continuity of $F$ everywhere, then the groupoid $\Omega $ must
be locally transitive.
\end{abstract}



\section{Introduction}

In order to establish notation, we give some definitions that can be found
(in equivalent form) in several places: \cite{ma},\cite{re1},\cite{mu}. To
define the term groupoid we begin with two sets, $\Omega $ and $B$, called
respectively the \emph{groupoid} and \emph{base}, together with two maps $%
\alpha $ and $\beta $, called respectively the \emph{source} and \emph{%
target projections}, a map $\varepsilon :u\rightarrow \tilde{u}$ from $B$ to 
$\Omega $ called the \emph{object inclusion map}, and a partial
multiplication $\left( x,y\right) \rightarrow xy$ in $\Omega $ defined on
the set 
\begin{equation*}
\Omega \ast \Omega =\left\{ \left( x,y\right) \in \Omega \times \Omega
:\alpha \left( x\right) =\beta \left( y\right) \right\} \text{(the set of 
\emph{composable pairs}),}
\end{equation*}
all subject to the following conditions:

\begin{enumerate}
\item  $\alpha \left( xy\right) =\alpha \left( y\right) $ and $\beta \left(
xy\right) =\beta \left( x\right) $ for all $\left( x,y\right) \in \Omega
\ast \Omega $.

\item  $x\left( yz\right) =\left( xy\right) z$ for all $x,y,z\in \Omega $
such that $\alpha \left( x\right) =\beta \left( y\right) $ and $\alpha
\left( y\right) =\beta \left( z\right) $.

\item  $\alpha \left( \tilde{u}\right) =\beta \left( \tilde{u}\right) =u$
for all $u\in B$.

\item  $x\widetilde{\alpha \left( x\right) }=x$ and $\widetilde{\beta \left(
x\right) }x=x$ for all $x\in \Omega $.

\item  Each $x\in \Omega $ has a (two-sided) inverse $x^{-1}$ such that $%
\alpha \left( x^{-1}\right) =\beta \left( x\right) $, $\beta \left(
x^{-1}\right) =\alpha \left( x\right) $ and $x^{-1}x=\widetilde{\alpha
\left( x\right) }$, $xx^{-1}=\widetilde{\beta \left( x\right) \text{. }}$
\end{enumerate}

The elements of the base $B$ will usually be denoted by letters as $%
u,\,v,\,w $ while arbitrary elements of the groupoid $\Omega $ will be
denoted by $x,\,y,\,z$.

The fibres of the target and the source projections will be denoted $\Omega
^{u}=\beta ^{-1}\left( \left\{ u\right\} \right) $ and $\Omega _{v}=\alpha
^{-1}\left( \left\{ v\right\} \right) $, respectively. More generally, given
the subsets $U$, $V$ $\subset B$, we define $\Omega ^{U}=\beta ^{-1}\left(
U\right) $, $\Omega _{V}=\alpha ^{-1}\left( V\right) $ and $\Omega
_{V}^{U}=\beta ^{-1}\left( U\right) \cap \alpha ^{-1}\left( V\right) $. The
relation $u\sim v$ iff $\Omega _{v}^{u}\neq \emptyset $ is an equivalence
relation on $B$. Its equivalence classes are called \emph{orbits} and the
orbit of a unit $u$ is denoted $\left[ u\right] $. A groupoid is called 
\emph{transitive} iff it has a single orbit. The set $\Omega _{u}^{u}$,
obviously a group under the restriction of the partial multiplication in $%
\Omega $, is called \emph{isotropy group} at $u$.

A \emph{subgroupoid} of $\Omega $ is a pair of subsets $\Omega ^{\prime
}\subseteq \Omega $ , $B^{\prime }\subseteq B$, such that $\alpha \left(
\Omega ^{\prime }\right) \subseteq B^{\prime }$, $\beta \left( \Omega
^{\prime }\right) \subseteq B^{\prime }$, $\tilde{u}\in \Omega ^{\prime }$
for all $u\in B^{\prime }$ and $\Omega ^{\prime }$ is closed under the
partial multiplication and inversion in $\Omega $. The \emph{inner
subgroupoid} (\emph{the isotropy group bundle}) of $\Omega $ is the
subgroupoid $G\Omega =\underset{u}{\cup }\Omega _{u}^{u}$. The \emph{base
subgroupoid} of $\Omega $ is the subgroupoid $\tilde{B}=\left\{ \tilde{u}%
:\,u\in B\right\} $. An element $\tilde{u}$ of $\tilde{B}$ is called the 
\emph{unity} corresponding to $u$.

A function $F:\Omega \rightarrow \Omega ^{\prime }$, where $\Omega $ and $%
\Omega ^{\prime }$ are groupoids, is called \emph{morphism} if for all
composable pairs $\left( x,y\right) \in \Omega \ast \Omega $, $\left(
F\left( x\right) ,F\left( y\right) \right) $ is a composable pair in $\Omega
^{\prime }$ and $F\left( xy\right) =F\left( x\right) F\left( y\right) $. It
is useful to note that $F\left( x^{-1}\right) =F\left( x\right) ^{-1}$ for
all $x\in \Omega $ and that $F\left( \tilde{u}\right) $ belongs to the base
subgroupoid of $\Omega ^{\prime }$ for all $\tilde{u}$ in the base
subgroupoid of $\Omega $.

Groupoids are generalizations of groups, but there are many other structures
which fit naturally into the study of groupoids. We give some basic examples
that will be needed in the subsequent considerations.

\textit{Examples of groupoids:}

\begin{enumerate}
\item  \textbf{Groups}: A group $G$ is a groupoid with base $B=\left\{
e\right\} $ (the unit element) and $G\ast G=G\times G$.

\item  \textbf{Spaces}. A space $X$ may be regarded as a groupoid on itself
with $\alpha =\beta =id_{X}$ and every element is a unity.

\item  \textbf{Equivalence relations}. Let $\mathcal{E}\subset X\times X$ be
an equivalence relation on the set $X$. We give $\mathcal{E}$ the structure
of a groupoid on $X$ in the following way: $\alpha $ is the projection onto
the second factor of $X\times X$ and $\beta $ is the projection onto the
first factor; the object inclusion map is $x\rightarrow \left( x,x\right) $
and the partial multiplication is $\left( x,y\right) \left( y^{\prime
},z\right) =\left( x,z\right) $ defined if and only if $y=y^{\prime }$. The
inverse of $\left( x,y\right) $ is $\left( y,x\right) $. Two extreme cases
deserve to be single out. If $\mathcal{E}=X\times X$, then $\mathcal{E}$ is
called the trivial groupoid on $X$, while if $\mathcal{E}=diag\left(
X\right) $, then $\mathcal{E}$ is called the co-trivial groupoid on $X$ (and
may be identified with the groupoid in example 2).

\item  \textbf{Transformation groups}. Let $G$ be a group acting on a set $X$
such that for $x\in X$ and $g\in G$, $gx$ denotes the transform of $x$ by $g$%
. Let $e$ be the unit element of $G$. Then $G\times X$ becomes a groupoid on 
$X$ in the following way: $\alpha $ is the projection onto the second factor
of $G\times X$ and $\beta $ is the action $G\times X$ $\rightarrow X$
itself; the object inclusion map is $x\rightarrow \left( e,x\right) $ and
the partial multiplication is $\left( g,x\right) \left( g^{\prime },y\right)
=\left( gg^{\prime },y\right) $ defined if and only if $x=g^{\prime }y$. The
inverse of $\left( x,g\right) $ is $\left( g^{-1},gx\right) $.
\end{enumerate}

A \emph{topological groupoid} consists of a groupoid $\Omega $ on $B$
together with topologies on $\Omega $ and $B$ such that the maps which
define the groupoid structure are continuous: the projections $\alpha $,$%
\beta :\Omega \rightarrow B$, the object inclusion map $\varepsilon
:B\rightarrow \Omega $, $\varepsilon \left( u\right) =\tilde{u}$, the
inversion map $x\rightarrow x^{-1}$ from $\Omega $ to $\Omega $, and the
partial multiplication from $\Omega \ast \Omega $ to $\Omega $, where $%
\Omega \ast \Omega $ has the subspace topology from $\Omega \times \Omega $,
are continuous.

As a consequence of the definition of a topological groupoid, the inversion $%
x\rightarrow x^{-1}$ is a homeomorphism, the object inclusion map is a
homeomorphism onto the base subgroupoid $\tilde{B}$, the projections are
identification maps. If $\Omega $ is Hausdorff, then $\tilde{B}$ is closed
in $\Omega $. If $B$ is Hausdorff, then $\Omega \ast \Omega $ is closed in $%
\Omega \times \Omega $.

Kirill Mackenzie raised in \cite{ma} (p. 31) the following question: given a
morphism $F:\Omega \rightarrow \Omega ^{\prime }$, where $\Omega $ and $%
\Omega ^{\prime }$ are topological groupoids and $F$ is continuous on a
neighborhood of \ the base in $\Omega $, is it true that $F$ is continuous
everywhere? \ Obviously, the answer is affirmative if the groupoid $\Omega $
is a group or if $\Omega $ is a co-trivial groupoid. Also if $\Omega $ is a
''principal ''topological groupoid, Mackenzie has proved that the answer is
affirmative (Proposition 1.21/p. 30 \cite{ma}) \ A \emph{''principal''
groupoid} (Definition 1.18/p. 28 \cite{ma}) is a topological \emph{transitive%
} groupoid $\Omega $ on $B$ which satisfies the following conditions:

\begin{enumerate}
\item  For each $v\in B$, the map $\beta _{v}:\Omega _{v}\rightarrow B$%
,\thinspace\ defined by $\beta _{v}\left( x\right) =\beta \left( x\right) $\
for all $x\in \Omega _{v}$, is open.

\item  For each $v\in B$, the map $\delta _{v}:\Omega _{v}\times \Omega
_{v}\rightarrow \Omega $,\thinspace\ defined by $\delta _{v}\left(
x,y\right) =xy^{-1}$\ for all $\left( x,y\right) \in \Omega _{v}\times
\Omega _{v}$, is open.
\end{enumerate}

We shall prove that the result remains true for a locally transitive
groupoid $\Omega $. By a locally transitive groupoid we shall mean a
topological groupoid $\Omega $ on $B$ for which the map $\beta _{v}:\Omega
_{v}\rightarrow B$ is open,\thinspace\ for each $v\in B$. Also we shall
prove that in general the answer to the question raised by Mackenzie is
negative. Moreover, we shall show that for a locally compact groupoid $%
\Omega $ with non-singleton orbits and having open target projection, if we
assume that the continuity of every morphism $F$ at all $\tilde{u}\in \tilde{%
B}$ implies the continuity of $F$ everywhere, then the groupoid $\Omega $
must be locally transitive. Furthermore, if we suppose that $\tilde{B}$ is
an open subset of $\Omega $ and the answer to the question of Mackenzie is
affirmative, then $B$ must be a discrete space. In particular, if $\Omega $
is the transformation group groupoid associated with the action of a
discrete group on a locally compact second countable Hausdorff space $B$,
then $\Omega $ has the property that every morphism which is continuous on a
neighborhood of the base must be globally continuous if and only if $B$ is a
discrete space.

\section{An example of a groupoid morphism which is continuous on a
neighborhood of the base but which is not continuous everywhere}

We shall construct a morphism (in the algebraic sense) $F:\Omega \rightarrow
\Omega ^{\prime }$ of topological groupoids with the property that the
continuity of $F$ on a neighborhood of the base in $\Omega $ does not imply
the continuity of $F$ everywhere on $\Omega $. Thus we give a negative
answer to a question raised by K. Mackenzie in \cite{ma} (p. 31).

Let $\mathbf{R}$ be the set of real numbers and $\mathbf{Z}$ the set of
integer numbers. Let $\Omega $ be the groupoid on $\mathbf{R}$ associated
with the following equivalence relation 
\begin{equation*}
\left\{ \left( x,y\right) \in \mathbf{R}^{2}:\;x-y\in \mathbf{Z}\right\}
\subset \mathbf{R}\times \mathbf{R}\text{.}
\end{equation*}

If $\Omega $ is endowed with the subspace topology from $\mathbf{R}\times 
\mathbf{R}$, $\Omega $ is a locally compact Hausdorff space. The base
subgroupoid of $\Omega $ 
\begin{equation*}
\tilde{B}=\left\{ \left( x,x\right) :\,x\in \mathbf{R}\right\}
\end{equation*}
is an open subset in $\Omega $. Let $\Omega ^{\prime }$ be the group of real
numbers under addition. Let $\phi :\mathbf{R}\rightarrow \mathbf{R}$ be
defined by $\phi \left( x\right) =0$ for all $x\leq \frac{1}{2}$ and $\phi
\left( x\right) =1$ otherwise. Let $F:\Omega \rightarrow \mathbf{R}$ be
defined by 
\begin{equation*}
F\left( x,y\right) =\phi \left( y\right) -\phi \left( x\right) \text{.}
\end{equation*}

For all $\left( \left( x,y\right) ,\left( y^{\prime },z\right) \right) \in
\Omega \ast \Omega $ we have 
\begin{eqnarray*}
F\left( \left( x,y\right) \left( y,z\right) \right) &=&F\left( x,z\right)
=\phi \left( z\right) -\phi \left( x\right) \\
&=&\phi \left( z\right) -\phi \left( y\right) +\phi \left( y\right) -\phi
\left( x\right) \\
&=&F\left( y,z\right) +F\left( x,y\right)
\end{eqnarray*}
Therefore $F$ is a morphism. Obviously, $F$ is continuous at every $\left(
x,x\right) $ in $\Omega $, but it is not continuous at $\left( \frac{1}{2},%
\frac{3}{2}\right) $.

\begin{remark}
The groupoid $\Omega $ is isomorphic and homeomorphic to the transformation
group groupoid obtained by letting $\mathbf{Z}$ act on $\mathbf{R}$. In this
setting the morphism (cocycle) $F$ is a coboundary determined by the
discontinuous function $\phi $.
\end{remark}

\section{Sufficient conditions for the continuity of a groupoid morphism}

\begin{theorem}
Let $\Omega $ be a topological groupoid on $B.$ Let us assume that for each $%
v\in B$, the map $\beta _{v}:\Omega _{v}\rightarrow B$,\thinspace\ defined
by $\beta _{v}\left( x\right) =\beta \left( x\right) $\ for all $x\in \Omega
_{v}$, is open. Let $\Omega ^{\prime }$ be a topological groupoid and $%
F:\Omega \rightarrow \Omega ^{\prime }$ be a morphism (in the algebraic
sense) which is continuous at every $\tilde{u}\in \tilde{B}$. Then $F$ is
continuous everywhere on $\Omega $.
\end{theorem}

\begin{proof}
Let $\left( x_{i}\right) _{i}$ be a net converging to $x$ in $\Omega $. Let $%
v=\alpha \left( x\right) $ and let us note that $\left( \beta \left(
x_{i}\right) \right) _{i}$ converges to $\beta \left( x\right) =\beta
_{v}\left( x\right) $ in $B$. Since $\beta _{v}:\Omega _{v}\rightarrow \ B$
is open, we may pass to a subnet and assume that there exists a net $\left(
y_{i}\right) _{i}$ in $\Omega _{v}$ converging to $x$ such that $\beta
_{v}\left( y_{i}\right) =\beta \left( x_{i}\right) $ for all $i$. The net $%
\left( y_{i}x^{-1}\right) _{i}$ converges to $xx^{-1}=\widetilde{\beta
\left( x\right) }$. From the continuity of $F$ at $\widetilde{\beta \left(
x\right) }$ it follows that $\lim\nolimits_{i}F\left( y_{i}x^{-1}\right)
=F\left( \widetilde{\beta \left( x\right) }\right) $, or equivalently $%
\lim\nolimits_{i}F\left( y_{i}\right) =F\left( x\right) $. On the other
hand, the net $\left( y_{i}^{-1}x_{i}\right) _{i}$ converges to $x^{-1}x=%
\widetilde{\alpha \left( x\right) }$. Now from the continuity of $F$ at $%
\widetilde{\alpha \left( x\right) }$ it follows that $\lim\nolimits_{i}F%
\left( y_{i}^{-1}x_{i}\right) =F\left( \widetilde{a\left( x\right) }\right) $%
. Thus 
\begin{eqnarray*}
\lim\nolimits_{i}F\left( x_{i}\right) &=&\lim\nolimits_{i}F\left( \left(
y_{i}y_{i}^{-1}\right) x_{i}\right) =\lim\nolimits_{i}F\left( y_{i}\left(
y_{i}^{-1}x_{i}\right) \right) \\
&=&\lim\nolimits_{i}F\left( y_{i}\right) F\left( y_{i}^{-1}x_{i}\right)
=F\left( x\right) F\left( \widetilde{a\left( x\right) }\right) \\
&=&F\left( x\right) \text{. }
\end{eqnarray*}
\end{proof}

\section{Necessary conditions for the continuity of a groupoid morphism}

Let $\Omega $ on $B$ be a locally compact (Hausdorff) groupoid. This means
that $\Omega $ is a topological groupoid whose topology is locally compact
(Hausdorff). A subset $A$ of $\Omega $ is called $\alpha $-relatively
compact if $A\cap \alpha ^{-1}\left( K\right) $ is relatively compact for
any compact subset $K$ of $B$. Jean Renault proved that if $\Omega $ is a
locally compact groupoid whose base $B$ is paracompact, then the base
subgroupoid $\tilde{B}$ admits a fundamental system of $\alpha $-relatively
compact neighborhoods (Proposition II.1.9/p. 56 \cite{re1}).

Throughout this section we shall fix a system of positive Radon measures
indexed by $B\times B$%
\begin{equation*}
\left\{ \lambda _{v}^{u},\,\left( u,v\right) \in B\times B\right\}
\end{equation*}
satisfying the following conditions:

\begin{enumerate}
\item  $supp\left( \lambda _{v}^{u}\right) =\Omega _{v}^{u}$ for all $u\sim
v $.

\item  $\lambda _{v}^{u}=0$ if $\Omega _{v}^{u}=\emptyset $.

\item  $\sup\limits_{u,v}\lambda _{v}^{u}\left( K\right) <\infty $ for all
compact $K\subset \Omega $.

\item  $\int f\left( y\right) d\lambda _{v}^{\beta \left( x\right) }\left(
y\right) $ =$\int f\left( xy\right) d\lambda _{v}^{\alpha \left( x\right)
}\left( y\right) $ for all $x\in \Omega $ and $v\,\in \left[ \,\beta \left(
x\right) \right] $.
\end{enumerate}

A construction of such system of measures can be found in \cite{re2}. In
Section $1$ of \cite{re2} Jean Renault constructs a Borel Haar system $%
\left\{ \lambda _{u}^{u},\,u\in B\right\} $ for the inner subgroupoid $%
G\Omega $. One way to do this is to choose a function $F_{0}$ continuous
with $\alpha $-relatively compact support which is nonnegative and equal to $%
1$ at each $\tilde{u}\in \tilde{B}.$ Then for each $u\in B$ choose a left
Haar measure $\lambda _{u}^{u}$ on $\Omega _{u}^{u}$ so the integral of $%
F_{0}$ with respect to $\lambda _{u}^{u}$ is $1.$

Renault defines $\lambda _{v}^{u}=x\lambda _{v}^{v}$ if $x\in \Omega
_{v}^{u} $ (where $x\lambda _{v}^{v}\left( f\right) =\int f\left( xy\right)
d\lambda _{v}^{v}\left( y\right) $, as usual). If $z$ is another element in $%
\Omega _{v}^{u}$, then $x^{-1}z\in \Omega _{v}^{v}$, and since $\lambda
_{v}^{v}$ is a left Haar measure on $\Omega _{v}^{v}$, it follows that $%
\lambda _{v}^{u}$ is independent of the choice of $x$. If $K$ is a compact
subset of $\Omega $, then $\sup\limits_{u,v}\lambda _{v}^{u}\left( K\right)
<\infty $.

\begin{lemma}
\label{lW}Let $\Omega $ be a locally compact groupoid on $B$, whose base
subgroupoid $\tilde{B}$ has a fundamental system of $\alpha $-relatively
compact neighborhoods. Let $v\in B$ and $f:G\rightarrow \ \mathbf{R}$ be a
continuous function with compact support. Then for each $\varepsilon >0$
there is a symmetric $a$-relatively compact neighborhood $W_{\varepsilon }$\
of $\tilde{B}$, such that: 
\begin{equation*}
\left| \int f\left( y\right) d\lambda _{v}^{\alpha \left( x\right) }\left(
y\right) \,-\int f\left( y\right) d\lambda _{v}^{\beta \left( x\right)
}\left( y\right) \right| <\varepsilon \lambda _{v}^{\alpha \left( x\right)
}\left( UL\right) \;\text{for all }\,x\in W_{\varepsilon }
\end{equation*}
where $L$ is the support of $f$ and $U$ is a symmetric $\alpha $-compact
neighborhood of $\tilde{B}$.
\end{lemma}

\begin{proof}
Since $\Omega $ is a locally compact groupoid whose base subgroupoid $\tilde{%
B}$ has a fundamental system of $\alpha $-relatively compact neighborhoods,
by Proposition 3.4/p. 40 \cite{bu}, it results that $f$ is left ''uniformly
continuous''. Hence for each $\varepsilon >0$ there is a symmetric $\alpha $%
-relatively compact neighborhood $W_{\varepsilon }$\ of $\tilde{B}$, such
that 
\begin{equation*}
x\in W_{\varepsilon }\,\,\;\Rightarrow \;\left| f\left( xy\right) \,-f\left(
y\right) \right| <\varepsilon \;\text{for all }\,y\in \Omega ^{\alpha \left(
x\right) }
\end{equation*}

Let $U$ be a symmetric $\alpha $-compact neighborhood of the base
subgroupoid $\tilde{B}$ and assume, without loss of generality, that $%
W_{\varepsilon }\subset U$. Then 
\begin{equation*}
\left| f\left( xy\right) -f\left( y\right) \right| <\varepsilon 1_{UL}\left(
y\right) \text{,}
\end{equation*}
for all $\,x\in W_{\varepsilon },\,\,y\in \Omega ^{\alpha \left( x\right) }$%
, where $1_{UL}$ denotes the indicator function of $UL$. It follows that 
\begin{equation*}
\left| \int f\left( xy\right) d\lambda _{v}^{\alpha \left( x\right) }\left(
y\right) -\int f\left( y\right) d\lambda _{v}^{\alpha \left( x\right)
}\left( y\right) \right| \leq \int \left| f\left( xy\right) -f\left(
y\right) \right| d\lambda _{v}^{\alpha \left( x\right) }\left( y\right)
\end{equation*}
\begin{equation*}
<\int \varepsilon 1_{UL}\left( y\right) d\lambda _{v}^{\alpha \left(
x\right) }\left( y\right) =\varepsilon \lambda _{v}^{\alpha \left( x\right)
}\left( UL\right) \;
\end{equation*}
for all $\,x\in W_{\varepsilon }$.
\end{proof}

\bigskip

\begin{proposition}
\label{cont}Let $\Omega $ be a locally compact groupoid on $B$, whose base
subgroupoid $\tilde{B}$ has a fundamental system of $\alpha $-relatively
compact neighborhoods. For each $v\in B$ and each continuous function with
compact support, $f:\Omega \rightarrow \ \mathbf{R}$, let us denote by $%
F_{f}^{v}:\Omega \rightarrow \ \mathbf{R}$ the map defined by 
\begin{equation*}
F_{f}^{v}\left( x\right) =\int f\left( y\right) d\lambda _{v}^{\alpha \left(
x\right) }\left( y\right) -\int f\left( y\right) d\lambda _{v}^{\beta \left(
x\right) }\left( y\right) \text{ for all }\,x\in \Omega \text{.}
\end{equation*}
Then

\begin{enumerate}
\item  $F_{f}^{v}$ is a morphism of groupoids.

\item  For each $\varepsilon >0$ there is a symmetric neighborhood $%
W_{\varepsilon }$ of the base subgroupoid $\tilde{B}$ such that: 
\begin{equation*}
\left| F_{f}^{v}\left( x\right) -F_{f}^{v}\left( \tilde{u}\right) \right|
=\;\left| F_{f}^{v}\left( x\right) \right| <\varepsilon \text{ for all }x\in
W_{\varepsilon }.
\end{equation*}

Consequently, $F_{f}^{v}$ is a morphism of groupoids which\ is continuous at
every unity $\tilde{u}\in \tilde{B}$.
\end{enumerate}
\end{proposition}

\begin{proof}
For all $\left( x,y\right) \in \Omega \ast \Omega $ we have 
\begin{eqnarray*}
F_{f}^{v}\left( xy\right) &=&\int f\left( z\right) d\lambda _{v}^{\alpha
\left( y\right) }\left( z\right) -\int f\left( z\right) d\lambda _{v}^{\beta
\left( x\right) }\left( z\right) \\
&=&\int f\left( z\right) d\lambda _{v}^{\alpha \left( y\right) }\left(
z\right) -\int f\left( z\right) d\lambda _{v}^{\beta \left( y\right) }\left(
z\right) + \\
&&\,\ \;\;\;\;\;\;\;\;\;+\int f\left( z\right) d\lambda _{v}^{\alpha \left(
x\right) }\left( z\right) -\int f\left( z\right) d\lambda _{v}^{\beta \left(
x\right) }\left( z\right) \\
&=&F_{f}^{v}\left( y\right) +F_{f}^{v}\left( x\right)
\end{eqnarray*}
Therefore $F_{f}^{v}$ is a morphism. Let $L$ be the support of $f$.

According to the preceding lemma, it follows that for each $\varepsilon >0$
there is a symmetric neighborhood $W_{\varepsilon }$ of the base subgroupoid 
$\tilde{B}$ such that: 
\begin{equation*}
\;\left| F_{f}^{v}\left( x\right) -F_{f}^{v}\left( \tilde{u}\right) \right|
<\varepsilon \lambda _{v}^{\alpha \left( x\right) }\left( UL\right) \text{
for all }x\in W_{\varepsilon }\,\,
\end{equation*}
\begin{equation*}
\;\left| F_{f}^{v}\left( x\right) -F_{f}^{v}\left( \tilde{u}\right) \right|
<\varepsilon \underset{u\in B}{\sup }\lambda _{v}^{u}\left( UL\right) \text{
for all }x\in W_{\varepsilon }
\end{equation*}
Hence $F_{f}^{v}$ is continuous at every $\tilde{u}$ $\in $ $\tilde{B}$.
\end{proof}

\begin{proposition}
\label{pE}Let $\Omega $ be a locally compact groupoid on $B$, whose base
subgroupoid $\tilde{B}$ has a fundamental system of $\alpha $-relatively
compact neighborhoods. Let $v\in B$. Then the following conditions are
equivalent:

\begin{enumerate}
\item  The map $\beta _{v}:\Omega _{v}\rightarrow B,\,\beta _{v}\left(
x\right) =\beta \left( x\right) $\ is open (or equivalently, the map $\alpha
_{v}:\Omega ^{v}\rightarrow B,\;a_{v}\left( x\right) =\alpha \left( x\right) 
$ is open).

\item  For each continuous function with compact support, $f:\Omega
\rightarrow \ \mathbf{R}$, the map 
\begin{equation*}
u\overset{\phi _{f}^{v}}{\rightarrow }\int f\left( y\right) d\lambda
_{v}^{u}\left( y\right) \;\left[ :B\rightarrow \mathbf{R}\right]
\end{equation*}
is continuous at every $u\in \left[ v\right] $.
\end{enumerate}
\end{proposition}

\begin{proof}
$1=>2$. Let $f:\Omega \rightarrow \ \mathbf{R}$ be a continuous function
with compact support. Let us denote by $F_{f}^{v}:\Omega \rightarrow \ 
\mathbf{R}$ the map defined by 
\begin{equation*}
F_{f}^{v}\left( x\right) =\int f\left( y\right) d\lambda _{v}^{\alpha \left(
x\right) }\left( y\right) -\int f\left( y\right) d\lambda _{v}^{\beta \left(
x\right) }\left( y\right) \text{ for all }\,x\in \Omega \text{.}
\end{equation*}
Let $u\in \left[ v\right] $ and let $\left( u_{i}\right) _{i}$ be a net in $B
$ converging to $u$. Let $x\in \Omega _{v}$ be such that $\beta _{v}\left(
x\right) =u$. Since $\beta _{v}:\Omega _{v}\rightarrow \ B$ is open, we may
pass to a subnet and assume that there exists a net $\left( x_{i}\right) _{i}
$ in $\Omega _{v}$ converging to $x$ such that $\beta _{v}\left(
x_{i}\right) =u_{i}$ for all $i$. Thus, we obtain a net $\left(
x_{i}x^{-1}\right) _{i}$ converging to $\tilde{u}$ in $\Omega $. \ According
to Proposition~\ref{cont}, $F_{f}^{v}$ is continuous at $\tilde{u}$. Hence 
\begin{equation*}
\lim\nolimits_{i}F_{f}^{v}\left( x_{i}x^{-1}\right) =0
\end{equation*}
that is 
\begin{equation*}
\lim\nolimits_{i}F_{f}^{v}\left( x_{i}\right) -F_{f}^{v}\left( x\right) =0%
\text{.}
\end{equation*}
Since $F_{f}^{v}\left( x_{i}\right) =\phi _{f}^{v}\left( \alpha \left(
x_{i}\right) \right) -\phi _{f}^{v}\left( \beta \left( x_{i}\right) \right)
=\phi _{f}^{v}\left( v\right) -\phi _{f}^{v}\left( u_{i}\right) $ and $%
F_{f}^{v}\left( x\right) =\phi _{f}^{v}\left( v\right) -\phi _{f}^{v}\left(
u\right) $, it follows that $\lim\nolimits_{i}\phi _{f}^{v}\left(
u_{i}\right) =\phi _{f}^{v}\left( u\right) $.

$2=>1$. Let $U$ be an open set in $\Omega $ and $u_{0}\in \beta _{v}\left(
U\cap \Omega _{v}\right) $. Let $x_{0}$ $\in $ $U$ be such that $u_{0}=\beta
\left( x_{0}\right) $. Let us choose a compactly supported nonnegative
continuous function, $f:\Omega \rightarrow \left[ 0,\,1\right] $, such that $%
f(x_{0})>0$ and $f$ vanishes outside $U$. Since the map 
\begin{equation*}
u\mapsto \int f\left( y\right) d\lambda _{v}^{u}\left( y\right) \;\left[
:B\rightarrow \mathbf{R}\right]
\end{equation*}
is continuous at $u_{0}$, it follows that 
\begin{equation*}
W=\left\{ u:\int f\left( y\right) d\lambda _{v}^{u}\left( y\right) >\frac{1}{%
2}\int f\left( y\right) d\lambda _{v}^{u_{0}}\left( y\right) \right\}
\end{equation*}
is a neighborhood of \ $u_{0}$. It is easy to see that $W$ is contained in
the set 
\begin{equation*}
\left\{ u:\int f\left( y\right) d\lambda _{v}^{u}\left( y\right) \neq
0\right\}
\end{equation*}
which is contained in $\beta _{v}\left( U\cap \Omega _{v}\right) $.
Therefore $\beta _{v}\left( U\cap \Omega _{v}\right) $ is a neighborhood of $%
u_{0}$.
\end{proof}

\begin{theorem}
\label{tF}Let $\Omega $ be a locally compact groupoid on $B$, whose base
subgroupoid $\tilde{B}$ has a fundamental system of $\alpha $-relatively
compact neighborhoods. Suppose that $\beta $, and hence $\alpha $, are open
maps. For each $v\in B$ and each continuous function with compact support, $%
f:\Omega \rightarrow \ \mathbf{R}$, let us denote by $F_{f}^{v}:\Omega
\rightarrow \ \mathbf{R}$ the map defined by 
\begin{equation*}
F_{f}^{v}\left( x\right) =\int f\left( y\right) d\lambda _{v}^{\alpha \left(
x\right) }\left( y\right) -\int f\left( y\right) d\lambda _{v}^{\beta \left(
x\right) }\left( y\right) ,\;\,x\in \Omega
\end{equation*}
and let us denote by $\phi _{f}^{v}$ : $B$ $\rightarrow $\ $\mathbf{R}$ the
map defined by 
\begin{equation*}
\phi _{f}^{v}\left( u\right) =\int f\left( y\right) d\lambda _{v}^{u}\left(
y\right) ,\;u\in B\text{.}
\end{equation*}
For each continuous function with compact support, $f:\Omega \rightarrow \ 
\mathbf{R}$, assume that $F_{f}^{v}$ is continuous on $\Omega $. If $\left[ v%
\right] $ is not singleton, then $\phi _{f}^{v}$ is continuous at every $%
u\in \left[ v\right] $.
\end{theorem}

\begin{proof}
Let $u\in B\cap \left[ v\right] $. Because $\left[ v\right] $ contains at
least two elements, there is $v_{0}\in \left[ v\right] $ such that $%
v_{0}\neq u$. Since $u\neq v_{0}$, there are two open subsets of $B$, $u\in
U_{0}$ and $v_{0}\in V_{0}$ such that $U_{0}$ $\cap $ $V_{0}$ $=$ $\emptyset 
$. Let $g:B\rightarrow \left[ 0,1\right] $ be a continuous function with the
following properties:

\begin{enumerate}
\item  $g\left( w\right) =1$ for all $w$ in a compact neighborhood $U$ of $u 
$ contained in $U_{0}$

\item  $g\left( w\right) =0$ for all $w\notin U_{0}$.
\end{enumerate}

For each continuous function with compact support, $f:\Omega _{v}\rightarrow
\ \mathbf{R}$, let us set $f_{U}$ $=$ $f$ $g\circ $ $\beta $. For each $w\in
B$, we have 
\begin{eqnarray*}
\phi _{f_{U}}^{v}\left( w\right) &=&\int f_{U}\left( y\right) d\lambda
_{v}^{w}\left( y\right) =\int f\left( y\right) g\left( \beta \left( y\right)
\right) d\lambda _{v}^{w}\left( y\right) = \\
&=&g\left( w\right) \int f\left( y\right) d\lambda _{v}^{w}\left( y\right)
=g\left( w\right) \phi _{f}^{v}\left( w\right)
\end{eqnarray*}
Therefore $\phi _{f}^{v}(w)=\phi _{f_{U}}^{v}\left( w\right) $ for all $w\in
U$ and $\phi _{f}^{v}(w)=0$ for all $w\in V_{0}$.

Let $\left( u_{i}\right) _{i}$ be a net converging to $u$. Without loss of
generality we may assume that $u_{i}\in U$ for all $i$, and hence, $\phi
_{f}\left( u_{i}\right) $ $=\phi _{f_{U}}\left( u_{i}\right) $ for all $i$.
So to prove that $\phi _{f}^{v}$ is continuous in $u$ it suffices to show
that $\lim_{i}\phi _{f_{U}}^{v}\left( u_{i}\right) $ $=$ $\phi
_{f_{U}}^{v}\left( u\right) $. Let us choose $x$ such that $\beta \left(
x\right) $ $=$ $u$ and $\alpha \left( x\right) $ $=$ $v_{0}$. Because $\beta
:\Omega \rightarrow \ B$ is open, we may pass to a subnet and assume that
there exists a net $\left( x_{i}\right) _{i}$ in $\Omega $ converging to $x$
such that $\beta \left( x_{i}\right) =u_{i}$ for all $i$. Since $%
F_{f_{U}}^{v}$ is continuous on $\Omega $, 
\begin{equation}
\ \lim\nolimits_{i}F_{f_{U}}^{v}\left( x_{i}\right) =F_{f_{U}}^{v}(x)=\phi
_{f_{U}}^{v}(v)-\phi _{f_{U}}^{v}(u)=-\phi _{f_{U}}^{v}(u)\text{.}
\label{e1}
\end{equation}

On the other hand, $\lim\nolimits_{i}\alpha \left( x_{i}\right) $ $=$ $%
\alpha \left( x\right) $ $=$ $v_{0}$ and for large $i$ we may assume that $%
\alpha \left( x_{i}\right) \in V_{0}$. Hence, $\phi _{f_{U}}^{v}(\alpha
\left( x_{i}\right) )=0$ and from $F_{f_{U}}^{v}\left( x_{i}\right) =\phi
_{f_{U}}^{v}(\alpha \left( x_{i}\right) )-\phi _{f_{U}}^{v}(\beta \left(
x_{i}\right) )$, it follows that 
\begin{equation}
F_{f_{U}}^{v}\left( x_{i}\right) =-\phi _{f_{U}}^{v}\left( \beta \left(
x_{i}\right) \right) \text{.}  \label{e2}
\end{equation}

By \ref{e1} and \ref{e2}, it results that $\lim\nolimits_{i}\phi
_{f_{U}}^{v}\left( u_{i}\right) =\phi _{f_{U}}^{v}\left( u\right) $, as
required.
\end{proof}

\bigskip

\begin{corollary}
Let $\Omega $ be a locally compact groupoid on $B$, whose base subgroupoid $%
\tilde{B}$ has a fundamental system of $\alpha $-relatively compact
neighborhoods. Suppose that $\beta $, and hence $\alpha $, are open maps.
For each $v\in B$ and each continuous function with compact support, $%
f:\Omega \rightarrow \ \mathbf{R}$, let us denote by $F_{f}^{v}:\Omega
\rightarrow \ \mathbf{R}$ the map defined by 
\begin{equation*}
F_{f}^{v}\left( x\right) =\int f\left( y\right) d\lambda _{v}^{\alpha \left(
x\right) }\left( y\right) -\int f\left( y\right) d\lambda _{v}^{\beta \left(
x\right) }\left( y\right) ,\;\,x\in \Omega \text{.}
\end{equation*}
Then each $F_{f}^{v}$ is a morphism of groupoids which\ is continuous at
every unity $\tilde{u}\in \tilde{B}$. If all morphisms $F_{f}^{v}$ are
continuous on $\Omega $, then for each $v\in B$ with non-singleton orbit,
the map $\beta _{v}:\Omega _{v}\rightarrow B,\,\beta _{v}\left( x\right)
=\beta \left( x\right) $\ is open.
\end{corollary}

\begin{proof}
According to the preceding theorem, for each $v$ with $\left[ v\right] \neq
\left\{ v\right\} $ and for each continuous function with compact support, $%
f:\Omega \rightarrow \ \mathbf{R}$, the map $\phi _{f}^{v}:B\rightarrow 
\mathbf{R}$ defined by 
\begin{equation*}
\phi _{f}^{v}\left( u\right) =\int f\left( y\right) d\lambda _{v}^{u}\left(
y\right) ,\;u\in B
\end{equation*}
is continuous at every $u\in \left[ v\right] $. Applying Proposition \ref{pE}%
, it follows that the map $\beta _{v}:\Omega _{v}\rightarrow B,\,\beta
_{v}\left( x\right) =\beta \left( x\right) $\ is open for all $v\in B$ with $%
\left[ v\right] \neq \left\{ v\right\} $.
\end{proof}

\begin{corollary}
Let $\Omega $ be a locally compact groupoid on $B$, for which the target
projection $\beta $ is an open map and whose base subgroupoid $\tilde{B}$ is
paracompact. Let us assume that for any morphism $F:\Omega \rightarrow 
\mathbf{R}$ the continuity of $F$ at all $\tilde{u}\in \tilde{B}$ implies
the continuity of $F$ everywhere on $\Omega $. Then for each $v\in B$ with
non-singleton orbit, the map $\beta _{v}:\Omega _{v}\rightarrow B,\,\beta
_{v}\left( x\right) =\beta \left( x\right) $\ is open.
\end{corollary}

\begin{proof}
If the base subgroupoid $\tilde{B}$ is paracompact, then $\tilde{B}$ has a
fundamental system of $\alpha $-relatively compact neighborhoods (see the
proof of Proposition II.1.9/p. 56 \cite{re1}). Therefore the hypotheses of
the preceding corollary are satisfied.
\end{proof}

\begin{definition}
A locally compact groupoid is $\beta $-discrete if its base subgroupoid is
an open subset (Definition I.2.6/p. 18 \cite{re1}).
\end{definition}

The $\beta $-discrete groupoids are generalizations of transformation groups 
$G\times B\rightarrow B$, where the group $G$ is assumed to be discrete. Any
locally compact groupoid for which $\beta $ (and hence $\alpha $) is a local
homeomorphism is $\beta $-discrete. Such groupoids are called \emph{etale
groupoids} (\cite{mu} p. 46) or \emph{sheaf groupoids}( \cite{ku} p. 209).

\begin{corollary}
Let $\Omega $ be a locally compact $\beta $-discrete groupoid, for which the
target projection $\beta $ is an open map and whose base subgroupoid $\tilde{%
B}$ is paracompact. Let us assume that for any morphism $F:\Omega
\rightarrow \mathbf{R}$ the continuity of $F$ on a neighborhood of $\tilde{B}
$ implies the continuity of $F$ everywhere on $\Omega $. If for each $v\in B$
the orbit $\left[ v\right] $ contains at least two elements, then\thinspace $%
B$ is a discrete space.
\end{corollary}

\begin{proof}
Applying Lemma I.2.7/p. 18 \cite{re1}, it follows that \ for each $v\in B$, $%
\Omega _{v}$ is a discrete space. According to the preceding corollary, \
for each $v\in B$, the map $\beta _{v}:\Omega _{v}\rightarrow B,\,\beta
_{v}\left( x\right) =\beta \left( x\right) $\ is open. Since $\,\beta
_{v}\left( \tilde{v}\right) =v$, it follows that $\left\{ v\right\} $ is
open in $B$.
\end{proof}

\begin{remark}
We have proved that for general topological groupoids the continuity of a
morphism on a neighborhood of the base does not necessarily imply the
continuity everywhere. Perhaps, the question of Mackenzie should be
reformulate: For what kind of groupoid morphism does the continuity on a
neighborhood of the base imply the continuity everywhere?
\end{remark}

\bigskip

University Constantin Br\^{a}ncu\c{s}i of T\^{a}rgu-Jiu

Bulevardul Republicii, Nr. 1, 210152 T\^{a}rgu-Jiu

Rom\^{a}nia

e-mail: ada@utgjiu.ro


\begin{thebibliography}{9}
\bibitem{bu}  M. Buneci, \textit{Haar systems and homomorphisms on groupoids}%
, Operator Algebras and Mathematical Physics: Conference Proceedings Constan%
\c{t}a (Romania), July 2-7, 2001, Theta Foundation (2003), 35-50.

\bibitem{ku}  A. Kumjian, \textit{On equivariant sheaf cohomology and
elementary }$C^{\ast }$\textit{-bundles}, J. Operator Theory, \textbf{20}
(1988), 207-240.

\bibitem{ma}  K. Mackenzie, \textit{Lie groupoids and Lie algebroids in
differential geometry}, LMS Lecture Note Ser., vol. \textbf{124}, Cambridge
Univ. Press, 1987.

\bibitem{mu}  P. Muhly, \textit{Coordinates in operator algebra}, (Book in
preparation).

\bibitem{re1}  J. Renault, \textit{A groupoid approach to }$C^{\ast }$%
\textit{- algebras}, Lecture Notes in Math. Springer-Verlag, \textbf{793}%
(1980).

\bibitem{re2}  J. Renault, \textit{The ideal structure of groupoid crossed
product }$C^{\ast }$\textit{- algebras}, J. Operator Theory, \textbf{25}
(1991), 3-36.
\end{thebibliography}
\end{document}